# Decomposition of the Catalan number into the sum of squares

Gennady Eremin
ergenns@gmail.com

June 11, 2020

**Abstract**. Analysis of the dynamics of the Dyck words helped solve the problem of representing the Catalan number as a sum of squares of natural numbers. In this case, the Dyck triangle is considered in different coordinates. In the calculations, we use the Catalan convolution matrix.

## Introduction

The well-formed parentheses, *Dyck words*, are sufficiently known. A Dyck word is a balanced string of left and right parentheses. The system of related parentheses determines the *Dyck dynamics*. For the Dyck word, the number of left and right parentheses is the same (the *first rule* of Dyck dynamics). In the initial subword, the number of right brackets never exceeds the number of left ones (the *second rule*). For the Dyck word of semilength $n$ ($n$ left parentheses and $n$ right ones), the second rule of Dyck dynamics is equivalent to the condition for the position $r_i$ of the $i$th right parenthesis [2]:

$$2i \le r_i \le n + i, \ 1 \le i \le n.$$

The number of Dyck words of semilength $n$ is equal to the $n$th Catalan number (see OEIS A000108). There are no restrictions to the length of the bracket sequences, so we can talk about the infinity of the set of Dyck words. This paper continues [1]. We will decompose the Catalan number into the sum of squares of integers.

## Dyck triangle

Usually the Dyck word is connected with a lattice path of *upsteps* of (1, 1), ascending diagonal vectors, and *downsteps* of (1, –1), descending diagonal vectors. The upstep corresponds to the left parenthesis, and the downstep corresponds to the right one. A *Dyck path* stays above the ground level (original line) that joins its initial and terminal points (vertices). Dyck paths are drawn in the first quadrant, each path starts at the origin. In such a lattice, the positions of the parentheses are plotted along the ground level, *i-axis*, and the *unbalance* (the excess number of left parentheses over right ones) is along the *j-axis*. The unbalance $j$ in each $i$th position can vary from 0 (the first $i$ brackets form a Dyck word) to $i$ (the first $i$ characters are left parentheses).

In [1] (see page 6), the author shows the path for the word ((()(()()))) of semi-length 6 (the 6-range). The polyline reaches a height of 4. It's not hard to see, for any achievable node $(i, j)$, the sum $i + j$ is an even number. Let us call an *n-triangle* figure with boundary points $(0, 0)$, $(n, n)$ and $(2n, 0)$. A Dyck path of semilength $n$ does not go beyond the $n$-triangle. The number of such paths is equal to the $n$th Catalan number.

For the node $(i, j)$, label $d(i, j)$ is equal to the number of paths from $(0, 0)$ to $(i, j)$. Inaccessible nodes are marked by 0. It is easy to see that $d(i, i) = 1, i \geq 0$, because there is only one path to $(i, i)$ (such paths correspond to adjacent left parentheses in the initial fragment of the Dyck word). For any node, its label (hereinafter *dynamics*) is determined by the *dynamics equation*:

$$d(i,j) = d(i-1, j-1) + d(i-1, j+1), \quad i \geq j > 0. \tag{1}$$

Using (1) and taking $d(0, 0) = 1$, it is easy to calculate the dynamics of all nodes. For example, you can only get to point $(2i, 0)$ from $(2i-1, 1)$, so

$$d(2i-1, 1) = d(2i, 0) = C_i, \quad i \leq n.$$

When Dyck paths are drawn, the abscissas of the diagonal vectors are indicated along the horizontal line. In the $i$th position of the Dyck word, the imbalance of the brackets does not exceed $i$, so on the plane, the set of non-zero labels forms the *Dyck triangle* (see Figure 1). Let's reformat the regular $n$-triangle.

| | | | | | | | | | | | | | | |
|---|---|---|---|---|---|---|---|---|---|---|---|---|---|---|
| 8 | | | | | | | | 1 | | 9 | | 54 | | 273 |
| 7 | | | | | | | 1 | | 8 | | 44 | | 208 | |
| 6 | | | | | | 1 | | 7 | | 35 | | 154 | | 637 |
| 5 | | | | | 1 | | 6 | | 27 | | 110 | | 429 | |
| 4 | | | | 1 | | 5 | | 20 | | 75 | | 275 | | 1001 |
| 3 | | | 1 | | 4 | | 14 | | 48 | | 165 | | 572 | |
| 2 | | 1 | | 3 | | 9 | | 28 | | 90 | | 297 | | 1001 |
| 1 | 1 | | 2 | | 5 | | 14 | | 42 | | 132 | | 429 | |
| 0 | 1 | 1 | | 2 | | 5 | | 14 | | 42 | | 132 | | 429 |
| | 0 | 1 | 2 | 3 | 4 | 5 | 6 | 7 | 8 | 9 | 10 | 11 | 12 | 13 | 14 |

Figure 1: The Dyck triangle.

In Figure 2 we showed a new $n$-triangle, an initial fragment of the Dyck triangle for $n = 6$. It is easy to see that for any node $x = (i, j)$, there is a *symmetric point* $y = (2n-i, j)$ with different dynamics in the General case. Further a symmetrical pair (*mirror*) $x, y$, we will denote $x \sim y$. In the center of the $n$-triangle, the $n$-column contains *self-symmetric nodes* $(n, j)$. In Figure 2, the 6th column contains four self-symmetric points 1, 5, 9, 5. It is easy to show that the sum of squares of these numbers is equal to $C_6$ (indeed, $1^2 + 5^2 + 9^2 + 5^2 = 132$).

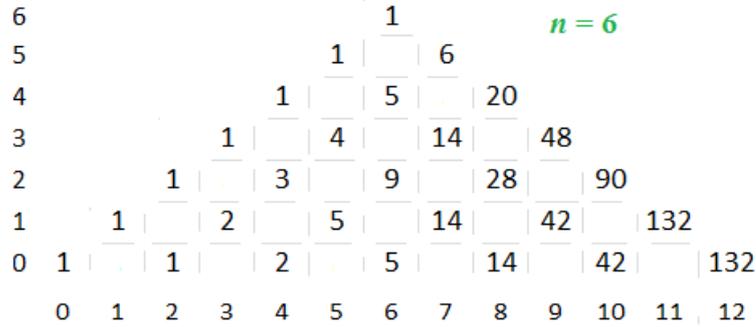

Figure 2: The 6-triangle.

In the *n*-triangle, let's invert the Dyck paths by placing their origin at $(2n, 0)$, and recount the *reverse dynamics* $\bar{d}(i, j)$. Then in the mirror $x \sim y$ the nodes are exchanged by the dynamics, that is, $\bar{d}(x) = d(y)$ and $\bar{d}(y) = d(x)$. For self-symmetric node the straight and reverse labels are the same, that is, $d(n, j) = \bar{d}(n, j)$, $j \leq n$.

Obviously, the number of Dyck paths through a node $x$ is equal to $d(x)\bar{d}(x)$. If we sum of such products for all nodes of any column in the *n*-triangle, we get the *n*th Catalan number. This implies the possibility of representing the Catalan number as the sum of squares of integers.

**Theorem 1.** *The nth Catalan number is equal to the sum of squares of elements in the nth column of the Dyck triangle.*

Let's call the squares of such a sum *Dyck squares*. It is easy to see that the number of Dyck squares in the sum is $\lfloor n/2 \rfloor + 1$ (the number of elements in the *n*th column), then

$$C_n = d^2(n, n) + d^2(n, n-2) + \ldots + d^2(n, n - 2\lfloor n/2 \rfloor).$$

Let's introduce the notation $t_{nk}$ ($= t_{n,k}$) for the elements of the *n*th column of the Dyck triangle:

$$C_n = \sum_{k=0}^{\lfloor n/2 \rfloor} t_{n,k}^2, \quad t_{n,k} = d(n, n - 2k). \qquad (2)$$

In (2) some terms are easy to determine, let's write these:

$t_{n, 0} = d(n, n) = 1$;
$t_{n, 1} = d(n, n-2) = n - 1$;
$t_{n, 2} = d(n, n-4) = 2 + 3 + \ldots + (n-3) + (n-2) = n(n-3)/2$;
$t_{n, \lfloor n/2 \rfloor} = d(n, n - 2\lfloor n/2 \rfloor) = C_{\lceil n/2 \rceil}$.

**Example 1**. Let's decompose $C_7 = 439$. In the sum, there are $\lfloor 7/2 \rfloor + 1 = 4$ terms that are calculated as follow:

$t_{7,0} = d(7, 7) = 1;$
$t_{7,1} = d(7, 5) = 7 - 1 = 6;$
$t_{7,2} = d(7, 3) = 7 \times (7-3)/2 = 14;$
$t_{7,3} = d(7, 1) = C_4 = 14.$

As a result, we get: $C_7 = 1^2 + 6^2 + 14^2 + 14^2 = 429.$ □

In the next section, we will try to obtain a General formula for $t_{nk}$.

## Convolution of Dyck triangle

Let us return to the infinite Dyck triangle (Figure 1). On the horizontal $i$-axis, the positions of Dyck words are indicated, i.e. an abscissa is the ordinal number of the current symbol (left or right parenthesis) of the word. On the vertical $j$-axis, the unbalance of brackets is indicated, i.e. an ordinate is the excess the number of left parentheses over right ones at the corresponding point on the $i$-axis. For a Dyck path of semilength $n$, the $i$-coordinate takes values from 0 to $2n$, and the $j$-coordinate takes values from 0 to $n$. For an arbitrary node $(i, j)$ of the Dyck triangle, the condition $j \leq i$ is mandatory and, in addition, $i + j =$ even.

In the Dyck triangle, nodes with the same abscissa are grouped as a column (a vertical line). There are $\lfloor i/2 \rfloor + 1$ elements in the $i$th column. Nodes with the same ordinate are grouped as an infinite row (a horizontal line). But there is a third way to select similar points; this is to group nodes in a diagonal.

It is easy to see that the nodes $(i, j)$, $i + j = 2n$, are located on the same descending diagonal line ($n$-isoline) with the upper point $(n, n)$ and the lower one $(2n, 0)$. The $n$-isoline is associated with the $n$th Catalan number, since the two bottom labels are equal to $C_n$.

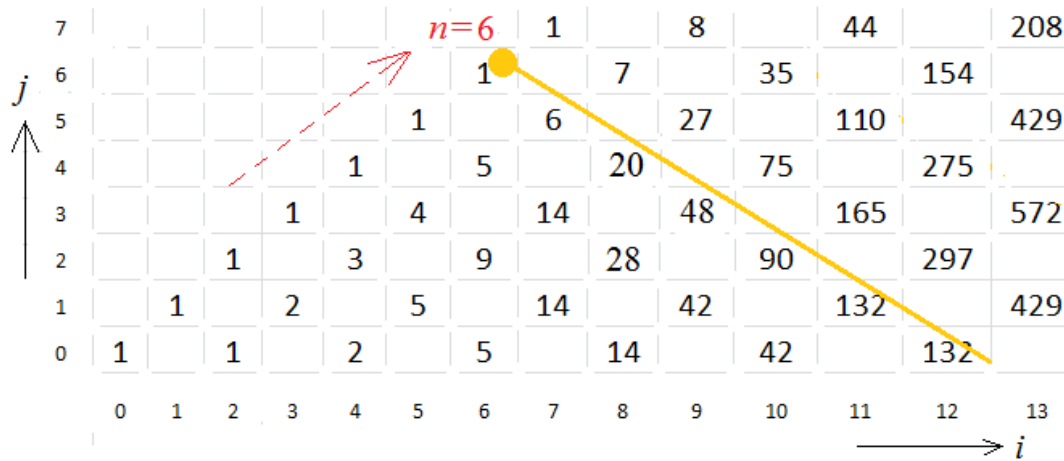

Figure 3: Turning the 6th diagonal isoline.

In the Figure 3 the 6-isoline is drawn in yellow and "fixed" with the button at the top. Again, in the Dyck triangle, each $n$-isoline is associated with the $n$th Catalan number. The converse is also true: for each natural number $n$, we can assign the $n$-isoline. This means that in addition to the two variables $i$ and $j$, we can use the third variable $n$, the index of the Catalan number. In Figure 3, the additional conditional (virtual) axis $n$ is shown with a red dashed line. All three variables are connected by the coordinate equation

$$i + j - 2n = 0, \ 0 \leq j \leq i \leq 2n. \tag{3}$$

In (3), the variable $n$ is, as before, the number of pairs of brackets in the Dyck word, but the status of $n$ has grown. Now $n$ is a coordinate along with $i$ and $j$. The grid $\{i, j\}$ is used in Figure 3, so $i$ and $j$ are independent variables, while the variable $n$ is dependent, $n = (i+j)/2$. Now we can say more specifically that in the coordinate grid $\{i, j\}$ we have the *Dyck ij-triangle*. But the grid can be changed.

Let's transform the Dyck *ij*-triangle into the *Dyck nj-triangle*, in other words, we need to change the coordinate grid $\{i, j\}$ to the coordinate grid $\{n, j\}$. This is simply done; it is enough to rotate each $k$-isoline to the vertical position, fixing the upper node $(k, k)$. The result is a *convolution* of the original Dyck triangle (see Figure 4).

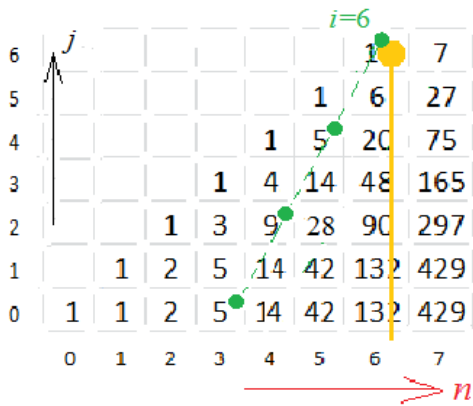
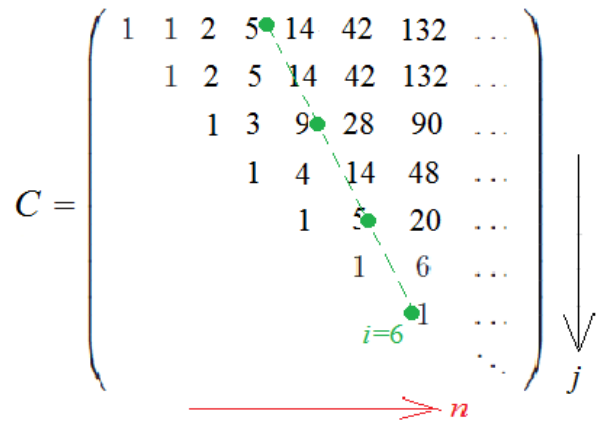

Figure 4.          Figure 5.

As we see the $i$-axis is replaced by the $n$-axis, and in our case the 6-isoline (from Figure 3) is converted in the 6th column. In the new coordinate grid, $n$ and $j$ are independent variables, while the variable $i$ is dependent, $i = 2n - j$. We continue to be interested in the nodes with the same coordinate $i$, as the labels of these points define Dyck squares in (2). So in Figure 4, the nodes with $i = 6$ are highlighted in green.

In the literature, you can find a numerical array similar to Figure 4. Figure 5 shows the known Catalan convolution matrix [4, 5]. The drawing is supplemented with axes $n, j$, and we highlighted the same four points with the coordinate $i = 6$. In the case $j = 0$ (upper row of the matrix), we get a set of Catalan numbers. It is easy to see that the Dyck *nj*-triangle and the Catalan convolution matrix are completely identical (the difference is only in the orientation of one axis).

Well-known formulas for Catalan numbers

$$C_n = \frac{(2n)!}{n!(n+1)!} = \frac{1}{n+1}\binom{2n}{n} = \binom{2n}{n} - \binom{2n}{n-1}, \ n \geq 0.$$

Elements of the Catalan convolution matrix (and also of the Dyck $nj$-triangle) are determined by the following formula [5, p. 2928]:

$$c(n, j) = \binom{2n-j}{n-j} - \binom{2n-j}{n-j-1}, \quad n, j \geq 0. \tag{4}$$

Let's define the relationship between the items of the Dyck $ij$-triangle and the Dyck $nj$-triangle. We go through all the highlighted nodes with $i = 6$ (see Figure 6):

$$t_{6,0} = d(6, 6) = c(6, 6) = 1; \quad t_{6,1} = d(6, 4) = c(5, 4) = 5;$$
$$t_{6,2} = d(6, 2) = c(4, 2) = 9; \quad t_{6,3} = d(6, 0) = c(3, 0) = 5.$$

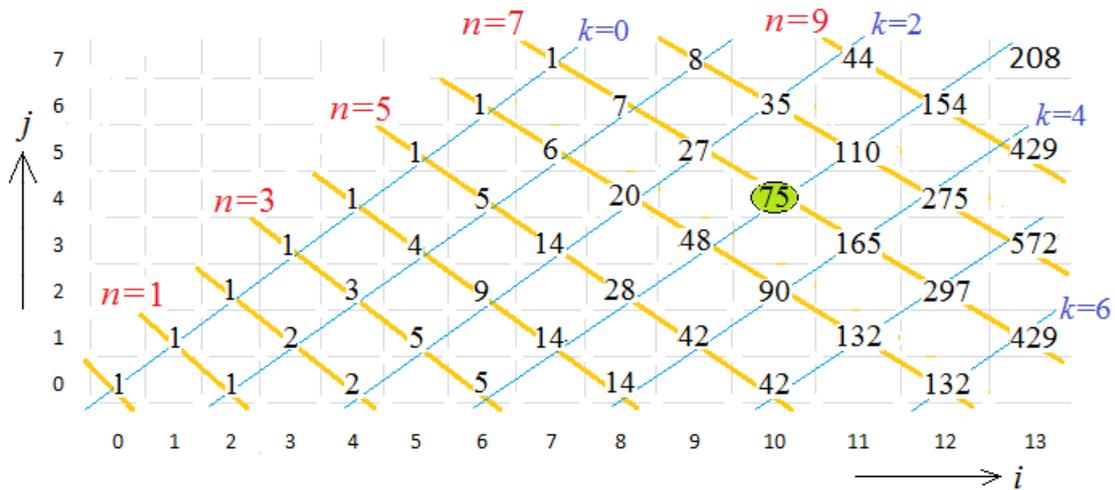

Figure 6: The Dyck $ij$-triangle with diagonal isolines.

It is not difficult to show that the relationship between the items of arrays is the following:

$$t_{ik} = d(i, j) = c(n, j), \quad n = (i+j)/2, \quad k = (i-j)/2. \tag{5}$$

To the variable $i, j, n$ we added the Dyck square index $k$. The inverse equalities are interesting:

$$i = n + k \quad \text{and} \quad j = n - k = i - 2k. \tag{6}$$

It is logical to consider (6) as an addition to (3). Figure 6 repeats the Dyck $ij$-triangle. Here the intersecting diagonals ($n = $ const and $k = $ const) give a visual representation of the functionality of the variables $n, k$. For example, the selected node $x = (10, 4)$ is at the intersection of the diagonals $n = 7$ and $k = 3$. The validity of (3) and (6) is obvious: $i_x = 7 + 3 = 10$, $j_x = 7 - 3 = 4$.

Making use of (4-6) and checking with Figure 6 (the green dashed line), we can write

$$t_{ik} = c(n, j) = c(i-k, i-2k) = \binom{2(i-k)-(i-2k)}{(i-k)-(i-2k)} - \binom{2(i-k)-(i-2k)}{(i-k)-(i-2k)-1}$$
$$= \binom{i}{k} - \binom{i}{k-1}.$$

Thus for the $n$th Catalan number we get the general formula of the Dyck square terms (2)

$$t_{nk} = \binom{n}{k} - \binom{n}{k-1} = n^{\underline{k-1}}(n-2k+1)/k!, \ 0 \le k \le \lfloor n/2 \rfloor, \qquad (7)$$

where $n^{\underline{k-1}} = n(n-1)\cdots(n-(k-2))$ is a falling factorial ($k-1$ factors).

The equality (7) is true for arbitrary $k$ [6, 7]. For example,

$t_{n,0} = n^{\underline{-1}} \times (n-0+1)/0! = \frac{1}{n+1} \times (n+1)/1 = 1,$
$t_{n,1} = n^{\underline{0}} \times (n-2+1)/1! = 1 \times (n-1)/1 = n-1,$
$t_{n,2} = n^{\underline{1}} \times (n-4+1)/2! = n \times (n-3)/2,$
$t_{n,3} = n^{\underline{2}} \times (n-6+1)/3! = n(n-1) \times (n-5)/6,$
$t_{n,4} = n^{\underline{3}} \times (n-8+1)/4! = n(n-1)(n-2) \times (n-7)/24,$ and so on.

In conclusion, we note that the Dyck $ij$-triangle and the Dyck $nj$-triangle (the Catalan convolution matrix) can be regarded as projections of some spatial construction, Dyck triangle in three-dimensional grid $\{i, j, n\}$ or *Dyck ijn-triangle*.

Gzhel State University, Moscow, 140155, Russia
http://www.en.art-gzhel.ru/